\documentclass[11pt,twoside,a4paper]{amsart}
\usepackage[english]{babel}

\usepackage[utf8]{inputenc}

\usepackage{pstricks-add}
\usepackage{array}
\usepackage{amsfonts}
\usepackage{amssymb}
\usepackage{amsmath}
\usepackage{graphicx}
\usepackage{graphics}
\usepackage{amsthm}

\newtheorem{theo}{Theorem}[section]
\newtheorem{de}[theo]{Definition}
\newtheorem{prop}[theo]{Proposition}
\newtheorem{lem}[theo]{Lemma}

\normalsize
\title{Optimal embedding of Meyer sets into model sets}

\author{Jean-Baptiste Aujogue}


\begin{document}

\maketitle

\begin{abstract}
We give a constructive proof that a repetitive Meyer multiple set of $\mathbb{R}^d$ admits a smallest model multiple set containing it colorwise.
\end{abstract}

\vspace{0.5cm}

Meyer sets are objects of central importance in the mathematical theory of Quasicrystals developed in the last thirty years. A special kind of Meyer sets is the class of so-called model sets, which are point patterns for which a geometric picture is available, the latter highly desirable for the understanding of the point pattern as well as for the computation of its relevant quantities. It is known for a long time that a Meyer set always embeds into at least one model set, but it is also true that a huge collection of model sets will contains a common prescribed Meyer set. Hence a natural question is whether there exist a "smallest" model set containing a given Meyer set, and our aim in this note is to answer this question, in the special case where the Meyer set is repetitive, by the affirmative. The existential result is as follows:

\vspace{0.3cm}

\textbf{Theorem.} \textit{For any repetitive Meyer set $\mathsf{\Lambda}$ of $\mathbb{R}^d$ is associated a unique model set $\underline{\mathsf{\Lambda}}$ such that whenever $\Delta $ is a model set with $\mathsf{\Lambda} \subseteq \Delta $ then one has}
\vspace{0.2cm}
\begin{align*} \mathsf{\Lambda} \subseteq \underline{\mathsf{\Lambda} } \subseteq \Delta 
\end{align*}

\vspace{0.3cm}

We provide a constructive proof of this result (see Theorem \ref{theo.principal}) by explicitly determining the CPS and window involved in the construction of $\underline{\mathsf{\Lambda} }$. The proof will be given in the slightly more general formalism of Meyer  multiple set, as it for instance naturally appears in the setting of substitution point sets.

\section{Meyer sets and model sets of $\mathbb{R}^d$}

A subset $\mathsf{\Lambda}$ of $\mathbb{R}^d$ is called a \textbf{Delone set} if it consists of a uniformly discrete collection of points, that is, it admits a uniform separation distance between any two of its points, and if it is relatively dense, meaning that any vector of $\mathbb{R}^d$ fits at uniformly bounded distance to some point of $\mathsf{\Lambda}$. A Delone set is called a \textbf{Meyer set}, or is said to have the Meyer property, if there is a finite subset $F$ of $\mathbb{R}^d$ such that
\vspace{0.2cm}
\begin{align*} \mathsf{\Lambda} - \mathsf{\Lambda} \, \subseteq \, \mathsf{\Lambda} +F
\end{align*}

\vspace{0.2cm}
Finally we call a Meyer multiple set a finite collection $\mathsf{\Lambda}= (\mathsf{\Lambda}_i)_{i\in I}$ of Meyer sets where the support $S(\mathsf{\Lambda}):= \cup _{i\in I} \mathsf{\Lambda}_i$ is again a Meyer set. We will often call such quantity a point pattern along the text. There are several very interesting equivalent formulations of the Meyer property which can be found in \cite{Moo}. The \textbf{hull} of a Meyer multiple set $\mathsf{\Lambda}$ consists of the collection $\mathbb{X}_{\mathsf{\Lambda}}$ of all Meyer multiples sets which locally resemble $\mathsf{\Lambda}$:
\vspace{0.2cm}
\begin{align*} \mathbb{X }_{\mathsf{\Lambda}}:= \left\lbrace \Lambda  \subset \mathbb{R}^d \, \vert \; \forall \, R> 0 \, \exists \, t\in \mathbb{R}^d :\;  \Lambda _i \cap B_R= (\mathsf{\Lambda}_i-t)\cap B_R \; \; \forall \, i\in I  \right\rbrace 
\end{align*}

\vspace{0.2cm}
where $B_R$ is the Euclidean ball of radius $R$ (or Euclidean $R$-ball). This collection is closed under the shift action of $\mathbb{R}^d$ by $\Lambda .t:= (\Lambda_i-t)_{i\in I}$, and admits a natural topology, the so-called \textbf{local topology} for which it is a compact space with continuous action of $\mathbb{R}^d$ (\cite{ LeeMoo, Sch}). Hence the Meyer multiple set $\mathsf{\Lambda}$ gives rise to a dynamical system $(\mathbb{X}_{\mathsf{\Lambda}}, \mathbb{R}^d)$.\\

We call a Meyer multiple set $\mathsf{\Lambda} $ \textbf{repetitive} whenever the dynamical system $(\mathbb{X}_{\mathsf{\Lambda}}, \mathbb{R}^d)$ is minimal. There is an intrinsic characterization of repetitivity on $\mathsf{\Lambda} $ (\cite{LaPl}) which will be of no use in our setting. There are subsets of $\mathbb{X}_{\mathsf{\Lambda}}$ which are of special interest, the so-called \textbf{canonical transversals} 
\vspace{0.2cm}
\begin{align*} \Xi _i:= \left\lbrace \Lambda  \in \mathbb{X}_{\mathsf{\Lambda}} \, \vert \; 0\in \Lambda _i \right\rbrace \; \;  \text{for} \, i\in I, \quad \qquad \Xi := \bigcup _{i \in I} \Xi _i = \left\lbrace \Lambda  \in \mathbb{X}_{\mathsf{\Lambda}} \, \vert \; 0\in S(\Lambda) \right\rbrace 
\end{align*}

\vspace{0.2cm}

We shall moreover consider another topology on $\mathbb{X}_{\mathsf{\Lambda}}$ which we call \textbf{combinatorial topology}, obtained from a uniformity for which a basis is given by
\vspace{0.2cm}
\begin{align*} \mathcal{U}_R := \left\lbrace (\Lambda , \Lambda ') \in \mathbb{X}_{\mathsf{\Lambda}}\times \mathbb{X}_{\mathsf{\Lambda}} \, \vert \; \Lambda _i \cap B_R= \Lambda '_i\cap B_R \; \; \forall \, i\in I  \right\rbrace 
\end{align*}

\vspace{0.2cm}
For all $R>0$. This topology is always strictly finer than the usual local topology, and endowing it on $\mathbb{X}_{\mathsf{\Lambda}}$ makes each transversal above a compact-open subset. Openness of the transversals for this topology is quite clear from the definition of the uniformity, whereas compactness is guaranteed by a consequence of the Meyer property called \textbf{finite local complexity} \cite{Auj, Ke}. The Meyer multiple set $\mathsf{\Lambda}$ is repetitive if and only if the $\mathbb{R}^d$-action on $\mathbb{X}_{\mathsf{\Lambda}}$ endowed with this new topology is minimal. 
\\

We shall now focus on a particular class of point patterns of $\mathbb{R}^d$: Suppose that we are given some lattice inside a space of the form $H\times \mathbb{R}^d$. Then one can form a point pattern by projecting along $\mathbb{R}^d$ the lattice points whose projection along $H$ belongs to some (nice) domain $W$ of $H$. The resulting point pattern is a so-called \textbf{model set} of $\mathbb{R}^d$. More precisely, we consider a triple $(H, \Gamma, s_H)$ called a \textbf{cut $\&$ project scheme} (or shortly \textbf{CPS}), where $H$ is a locally compact $\sigma$-compact Abelian group (LCA group for short), $\Gamma$ a countable subgroup of $\mathbb{R}^d$ and
\vspace{0.2cm}
\begin{align*}\begin{psmatrix}[colsep=1.5cm,
rowsep=0cm]
\Gamma  \; \; & \; \; H \\
\psset{arrows=->,linewidth=0.2pt, labelsep=1.5pt
,nodesep=0pt}
\ncline{1,1}{1,2}^{s_H}
\end{psmatrix}\end{align*}

\vspace{0.2cm}

a group morphism with a range $ s_H(\Gamma)$ dense in $H$, and furthermore with a graph $\mathcal{G}(s_H):= \left\lbrace (s_H(\gamma),\gamma)\in H \times \mathbb{R}^d\, \vert \, \gamma\in \Gamma\right\rbrace$ that it a \textit{lattice}, that is, a discrete and co-compact subgroup of $H\times \mathbb{R}^d$. The LCA group $H$ is commonly called the \textbf{internal space}, the subgroup $\Gamma$ of $\mathbb{R}^d$ the \textbf{structure group} of the CPS, and the morphism $s_H$ the \textbf{*-map} of the CPS. A model set issued from this CPS is then a point pattern of $\mathbb{R}^d$ obtained as
\vspace{0.2cm}
\begin{align*}\mathfrak{P}_{H}(W):= \left\lbrace \gamma \in \Gamma \; \vert \; s_H(\gamma ) \in W \right\rbrace 
\end{align*}

\vspace{0.2cm}

where $W$ is a compact topologically regular subset of $H$, that is, a compact set which is the closure of its interior $\mathring{W}$ in $H$. A model multiple set is a finite collection $(\mathfrak{P}(W_i))_{i\in I}$ of model sets issuing from a common CPS but with possibly different sets $W_i$. The finite collection $( W_i ) _{i\in I}$ of compact topologically regular subsets used to form a model multiple set, or the set $W$ in case we deal with a single model set, is called a \textbf{window}. The formal definition of a model multiple set now comes as follows:

\vspace{0.2cm}

\begin{de}\label{defmodelset} A model multiple set $\mathsf{\Lambda}= (\mathsf{\Lambda}_i)_{i\in I}$ is a point pattern such that there exists a CPS $(H, \Gamma, s_H)$ and a window $\left\lbrace  W_i \right\rbrace _{i\in I}$ in $H$ such that
\vspace{0.2cm}
\begin{align*} \mathfrak{P}_{H}( \stackrel{\circ }{W_i}) \, \subseteq \, \mathsf{\Lambda}_i \, \subseteq  \,\mathfrak{P}_{H}(W_i)
\end{align*}
\end{de}
\vspace{0.2cm}
In case the window admits boundary sets $ \partial W_i$ of null Haar measure in $H$ the resulting model multiple sets are said to be \textit{regular}. We are not assuming this here. It is a standard fact (\cite{Moo2}) that a model multiple set is always a Meyer multiple set.

\vspace{0.2cm}

\section{A general construction}

Let us consider throughout this section a fixed Meyer multiple set $\mathsf{\Lambda}$ with hull $\mathbb{X}_{\mathsf{\Lambda}}$. One can then provide a "model set-like" description of $\mathsf{\Lambda}$ as follows:\\

Let $\mathsf{\Gamma}$ be the subgroup of $\mathbb{R}^d$ generated by the support of $\mathsf{\Lambda}$. The subset $S(\mathsf{\Lambda})$ is Meyer and thus \textbf{finitely generated} (see \cite{LaPl}), which in particular means that $\mathsf{\Gamma}$ is countable. Consider the collection of point patterns in the hull $\mathbb{X}_{\mathsf{\Lambda}}$ having support into $\mathsf{\Gamma}$,
\begin{align*} \Xi ^{\mathsf{\Gamma}} := \left\lbrace \Lambda \in \mathbb{X}_{\mathsf{\Lambda}} \, \vert \; S(\Lambda )\subseteq \mathsf{\Gamma} \right\rbrace 
\end{align*}

\vspace{0.2cm}

this collection equipped with combinatorial topology. The topological space $\Xi ^{\mathsf{\Gamma}}$ hence defined contains the transversal $\Xi$, and in particular all the colored transversals $\Xi_i$ with $i\in I$: For, if we are given any $\Lambda \in \Xi$ then $0\in S(\Lambda)$ and thus any $\gamma \in S(\Lambda)$ provides  a two-point pattern $\left\lbrace 0, \gamma\right\rbrace $ of $S(\Lambda)$. As any bounded pattern of $\Lambda$ appears somewhere in $\mathsf{\Lambda}$ there must be some vector $t\in \mathbb{R}^d$ with $\left\lbrace t, t+\gamma\right\rbrace \subset S(\mathsf{\Lambda})$. In particular $\gamma$ lies into the difference set $S(\mathsf{\Lambda})-S(\mathsf{\Lambda})$, in turns contained into $\mathsf{\Gamma} -\mathsf{\Gamma} = \mathsf{\Gamma}$. This provides $S(\Lambda) \subseteq \mathsf{\Gamma}$, which shows that $\Lambda \in \Xi ^{\mathsf{\Gamma}}$ whenever $\Lambda \in \Xi$, as desired.\\

The space $\Xi ^{\mathsf{\Gamma}}$ with combinatorial topology is locally compact $\sigma$-compact \cite{Auj}, and is in fact a totally disconnected space. It is clear that if $\Lambda \in  \Xi ^{\mathsf{\Gamma}}$ then $\Lambda -\gamma \in  \Xi ^{\mathsf{\Gamma}}$ whenever $\gamma \in \mathsf{\Gamma}$, providing a shift action of $\mathsf{\Gamma}$. It is not difficult to show that, since $\mathsf{\Lambda}$ is assumed to be repetitive, the $\mathsf{\Gamma}$-action on the space $\Xi ^{\mathsf{\Gamma}}$ with combinatorial topology is minimal. It is also clear that the point pattern $\mathsf{\Lambda}$ itself falls into $\Xi ^{\mathsf{\Gamma}}$, and all its $\mathsf{\Gamma}$-translates as well. Therefore we have a locally compact $\sigma$-compact space $\Xi ^{\mathsf{\Gamma}}$, the countable subgroup $\mathsf{\Gamma}$ of $\mathbb{R}^d$ and moreover the mapping
\vspace{0.2cm}
\begin{align*}\begin{psmatrix}[colsep=1.5cm,
rowsep=0cm]
\mathsf{\Gamma} \ni \gamma  \; \; & \; \;s_{\Xi ^{\mathsf{\Gamma}}}(\gamma):= \mathsf{\Lambda}-\gamma \in \Xi ^{\mathsf{\Gamma}} 
\psset{arrows=|->,linewidth=0.2pt, labelsep=1.5pt
,nodesep=0pt}
\ncline{1,1}{1,2}
\end{psmatrix}\end{align*}

\vspace{0.2cm}

We may interpret the triple $(\Xi ^{\mathsf{\Gamma}}, \mathsf{\Gamma}, s_{\Xi ^{\mathsf{\Gamma}}})$ as a "CPS" associated with $\mathsf{\Lambda}$, with $\Xi ^{\mathsf{\Gamma}}$ as "internal space", $\mathsf{\Gamma}$ as "structure group" and $s_{\Xi ^{\mathsf{\Gamma}}}$ as "$*$-map". In addition the transversals $\Xi _i$, $i\in I$, are compact and open subsets of $\Xi ^{\mathsf{\Gamma}}$ and thus topologically regular, giving a "window" $\left\lbrace \Xi _i\right\rbrace _{i\in I}$ in this space. The point pattern $\mathsf{\Lambda}$ now appears as a "model set" associated with $(\Xi ^{\mathsf{\Gamma}}, \mathsf{\Gamma}, s_{\Xi ^{\mathsf{\Gamma}}})$ and $\left\lbrace \Xi _i\right\rbrace _{i\in I}$ in the sense that for each $i\in I$ one has
\vspace{0.2cm}
\begin{align*} \mathsf{\Lambda}_i = \mathfrak{P}_{\Xi ^{\mathsf{\Gamma}}}(\Xi _i):= \left\lbrace \gamma \in \mathsf{\Gamma} \, \vert \, s_{\Xi ^{\mathsf{\Gamma}}}(\gamma ) \in \Xi _i\right\rbrace  
\end{align*}

\vspace{0.2cm}

It is in fact possible to derive a true CPS and window form what is given above, and the key ingredient here is a certain relation on $\mathbb{X}_{\mathsf{\Lambda}}$ called, after its introduction in \cite{BaKe}, the \textbf{strong regional proximality} relation: Two point patterns $\Lambda , \Lambda '\in \mathbb{X}_{\mathsf{\Lambda}}$ are strongly regionally proximal, briefly denoted $\Lambda \sim _{srp} \Lambda '$, if for each radius $R$ there are $\Lambda _1, \Lambda _2$ and $t\in \mathbb{R}^d$ such that colorwise
\vspace{0.2cm}
\begin{align*} \Lambda \cap B_R &=  \Lambda _1 \cap B_R \\ \Lambda ' \cap B_R &= \Lambda _2 \cap B_R \\ (\Lambda _1-t)\cap B_R &=  (\Lambda _2-t) \cap B_R
\end{align*}

\vspace{0.2cm}

It is true that for a repetitive Meyer multiple set $\mathsf{\Lambda}$ this relation is a closed $\mathbb{R}^d$-invariant equivalence relation on $\mathbb{X}_{\mathsf{\Lambda}}$, something by no means obvious. In this case this relation is related with spectral features, as two point patterns are strongly regionally proximal if and only if they are undistinguished by the continuous eigenfunction of the dynamical system $(\mathbb{X}_{\mathsf{\Lambda}}, \mathbb{R}^d)$, see \cite{BaKe} for details.

\vspace{0.2cm}

\begin{prop}\label{prop.construction.CPS}\cite{Auj}(CPS construction) Let $\mathsf{\Lambda}$ be a repetitive Meyer multiple set of $\mathbb{R}^d$, and $\mathsf{\Gamma}$ be the countable subgroup of $\mathbb{R}^d$ generated by its support $S(\mathsf{\Lambda})$. Then:

\vspace{0.3cm}

$(i)$ If $\Lambda \in \Xi ^{\mathsf{\Gamma}}$ then its strong regional proximality class $[\Lambda ]_{srp}$ is contained into $\Xi ^{\mathsf{\Gamma}}$.

\vspace{0.25cm}

$(ii)$ The quotient space $\mathsf{H}:= \Xi ^{\mathsf{\Gamma}}\diagup \sim _{srp}$ with quotient topology admits a LCA group structure, such that one has a group morphism
\begin{align*}\begin{psmatrix}[colsep=1.5cm,
rowsep=0cm]
\mathsf{\Gamma} \ni \gamma  \; \; & \; \;s_{\mathsf{H}}(\gamma):= \left[ \mathsf{\Lambda}+\gamma\right] _{srp}\in \mathsf{H}
\psset{arrows=|->,linewidth=0.2pt, labelsep=1.5pt
,nodesep=0pt}
\ncline{1,1}{1,2}
\end{psmatrix}\end{align*}
\end{prop}

\vspace{0.2cm}

The above proposition showed that the triple $(\Xi ^{\mathsf{\Gamma}}, \mathsf{\Gamma}, s_{\Xi ^{\mathsf{\Gamma}}})$ we considered earlier yields by modding out the strong regional proximality relation a new triple $(\mathsf{H}, \mathsf{\Gamma}, s_{\mathsf{H}})$, with $\mathsf{H}$ an LCA group, $\mathsf{\Gamma}$ a countable subgroup of $\mathbb{R}^d$ and $s_{\mathsf{H}}$ a group morphism from the latter into the former. The following theorem provides a CPS and window to the Meyer multiple set $\mathsf{\Lambda}$ which will be of central importance for our concern:

\vspace{0.3cm}

\begin{theo}\label{theo.construction.CPS}\cite{Auj} Let $\mathsf{\Lambda}$ be a repetitive Meyer multiple set of $\mathbb{R}^d$, and $\mathsf{\Gamma}$ be the countable subgroup of $\mathbb{R}^d$ generated by its support $S(\mathsf{\Lambda})$. Then:

\vspace{0.3cm}

$(i)$ The triple $(\mathsf{H}, \mathsf{\Gamma}, s_{\mathsf{H}})$ is a CPS.

\vspace{0.25cm}

$(ii)$ Each subset $\mathsf{W}_i:= -\left[ \Xi _i\right] _{srp}$ of $\mathsf{H}$ is compact and topologically regular in $\mathsf{H}$ and thus form a window $\left\lbrace \mathsf{W}_i \right\rbrace _{i \in I}$ in $\mathsf{H}$.

\vspace{0.25cm}

$(iii)$ One has for each $i\in I$ that $\displaystyle \mathsf{\Lambda}_i \, \subseteq \, \mathfrak{P}_{\mathsf{H}}(\mathsf{W}_{i})= \bigcup _{\Lambda \in [\mathsf{\Lambda}]_{srp}} \Lambda _i$.

\end{theo}

\vspace{0.3cm}



\section{Statement and proof of the theorem}

What we aim to prove in this note is the following:

\vspace{0.2cm}

\begin{theo}\label{theo.principal} Let $\mathsf{\Lambda}$ be a repetitive Meyer multiple set of $\mathbb{R}^d$, and consider the associated CPS $(\mathsf{H} ,\mathsf{\Gamma } , s_{\mathsf{H}})$ and window $(\mathsf{W}_i)_{i\in I}$ provided by Theorem \ref{theo.construction.CPS}. Then the point pattern colorwise given by
\vspace{0.2cm}
\begin{align*} \underline{\mathsf{\Lambda}}_i := \mathfrak{P}_{\mathsf{H}}(\mathring{\mathsf{W}}_{i})\cup \mathsf{\Lambda}_i
\end{align*}
is a model multiple set in the sense of Definition \ref{defmodelset}, and has the property that whenever $\Delta $ is a model multiple set with $\mathsf{\Lambda} \subseteq \Delta $ colorwise then one has colorwise
\vspace{0.2cm}
\begin{align*} \mathsf{\Lambda} \subseteq \underline{\mathsf{\Lambda} } \subseteq \Delta 
\end{align*}
\end{theo}

\vspace{0.3cm}

The fact that $\underline{\mathsf{\Lambda}}$ is a model multiple set follows from part $(iii)$ of Theorem \ref{theo.principal}, as it gives $ \mathfrak{P}_{\mathsf{H}}(\mathring{\mathsf{W}}_{i})\, \subseteq \, \underline{\mathsf{\Lambda}}_i \, \subseteq \, \mathfrak{P}_{\mathsf{H}}(\mathsf{W}_{i})$. The subsequent part of this note is devoted to the proof of the above theorem. Let us begin by supposing that $\Delta $ is a model multiple set with $\mathsf{\Lambda}\subseteq \Delta $ colorwise. Therefore there is a CPS $(H^1, \Gamma ^1, s_{H_1})$ and a window $\left( W_i^1 \right) _{i\in I}$ in $H^1$ such that according to Definition \ref{defmodelset}
\vspace{0.2cm}
\begin{align}\label{interpolation} \mathfrak{P}_{H^1}(\mathring{W}_i^1) \subseteq \Delta _i \subseteq \mathfrak{P}_{H^1}( W_i^1)
\end{align}

\vspace{0.3cm}

We now make the two following observations: First the group $\mathsf{\Gamma}$ generated by the support of $\mathsf{\Lambda}$ is contained into $\Gamma^1$, and second each set $s_{H^1}(\mathsf{\Lambda}_i)$ is contained in the compact set $W_i^1$ and thus has compact closure $V_i^1$ in $H^1$.\\

\textbf{Step 1: Factoring out the redundancies.} We shall now modify the CPS $( H^1, \Gamma^1, s_{H^1})$ by modding out the redundancy subgroup
\vspace{0.2cm}
\begin{align*} \mathfrak{R}:= \left\lbrace w\in H^1 \, \vert \, V_i^1 +w=V_i^1 \; \forall \, i\in I\right\rbrace 
\end{align*}

\vspace{0.2cm}

associated with the collection $(V_i^1)_{i\in I}$ of compact subsets in $H^1$. This provides a LCA group $H^2:= H^1\diagup \mathfrak{R}$ and a new CPS $(H^2 ,\Gamma ^2, s_{H^2})$ where the structure group $\Gamma ^2:= \Gamma ^1$ remains unchanged, with $*$-map $s_{H^2}$ given by
\vspace{0.2cm}
\begin{align}\label{morphism2}\begin{psmatrix}[colsep=1.5cm,
rowsep=1.2cm]
\Gamma^2 \ni \gamma \; \; & \; \; s_H^2(\gamma):= [s_{H^1}(\gamma)]_\mathfrak{R} \in \; H^2
\psset{arrows=|->,linewidth=0.2pt, labelsep=1.5pt
,nodesep=0pt}
\ncline{1,1}{1,2}
\end{psmatrix}\end{align}

\vspace{0.2cm}

The proof that it forms an actual CPS comes from the fact that $\mathfrak{R}$ is compact, see also the discussion of \cite{LeeMoo}, section $5$ therein. Each compact subset $ V_i^2:= [V_i^1]_\mathfrak{R}$ is equal to the closure of $s_{H^2}(\mathsf{\Lambda}_i)$ in $H^2$, and form an irredundant family $(V_i^2)_{i \in I}$ in $H^2$.\\

\textbf{Step 2: Shrinking the structure group.} The second operation we wish to apply is shrinking the structure group $\Gamma$ we are considering to the group $\mathsf{\Gamma}$ generated by the support of $\mathsf{\Lambda}$. Since the structure group $\Gamma^2$ contains $\mathsf{\Gamma}$ we shall consider the LCA group $H$ to be the closure of $s_{H^2}(\mathsf{\Gamma})$ in $H^2$. Restricting the group morphism $s_{H^2}$ on $\mathsf{\Gamma}$ yields a group morphism
\vspace{0.2cm}
\begin{align}\label{morphism}\begin{psmatrix}[colsep=1.5cm,
rowsep=1.2cm]
\mathsf{\Gamma} \; \; & \; \; \; H
\psset{arrows=->,linewidth=0.2pt, labelsep=1.5pt
,nodesep=0pt}
\ncline{1,1}{1,2}^{s_H}
\end{psmatrix}\end{align}

\vspace{0.2cm}

Moreover each set $V_i^2$, being the closure of $s_{H^2}(\mathsf{\Lambda}_i)$ in $H^2$, is contained into $H$ and we shall rewrite each as $V_i$. It naturally comes that each $V_i$ is compact and the closure of $s_H(\mathsf{\Lambda}_i)$ in $H$.

\vspace{0.2cm}

\begin{prop}\label{lem.small.CPS} The triple $(H, \mathsf{\Gamma}, s_H)$ forms a CPS.
\end{prop}

\vspace{0.2cm}

\begin{proof} By construction the group morphism $s_H$ has a dense range $s_H(\mathsf{\Gamma})$ in the LCA group $H$. Moreover, $s_H$ has a uniformly discrete graph $\mathcal{G}(s_H)$ in $H\times\mathbb{R}^d $ since $\mathcal{G}(s_H)\subseteq \mathcal{G}(s_{H^2})$ which is uniformly discrete in $H^2\times\mathbb{R}^d $. What thus remains to show is the relative density of $\mathcal{G}(s_H)$ in $H\times\mathbb{R}^d $. Let $i\in I$ be chosen arbitrarily. Then the closure $V_i$ of $s_H(\mathsf{\Lambda}_i)$ is compact in $H$ and therefore one has an open relatively compact subset $U$ of $H$ with $V_i\subseteq U$. The point pattern $\mathsf{\Lambda}_i$ is relatively dense in $\mathbb{R}^d$, with radius $R$ of relative density. Then the graph of $s_H$ obeys
\vspace{0.2cm}
\begin{align*} \mathcal{G}(s_H) + (U-U)\times B_R = H\times\mathbb{R}^d 
\end{align*}

\vspace{0.2cm}

where $B_R$ is the Euclidean $R$-ball and $(U-U)\times B_R$ relatively compact in $H\times\mathbb{R}^d$. For, given $(w,t)\in H\times\mathbb{R}^d $ one has by density of $s_H(\mathsf{\Gamma})$ in $H$ a $\gamma\in \mathsf{\Gamma}$ with $w-s_H(\gamma)\in U$. Since $\mathsf{\Lambda}_i$ is $R$-relatively dense there is a $\gamma' \in \mathsf{\Lambda}_i$ with $\gamma'\in (t-\gamma)+B_R$. As $\gamma'$ lies into $\mathsf{\Lambda}_i$ we have $s_H(\gamma')\in U$. Now let $\gamma _0:= \gamma + \gamma ' \in \mathsf{\Gamma} + \mathsf{\Lambda}_i = \mathsf{\Gamma}$: Then $w \in s_H(\gamma)+ U \subset s_H(\gamma)+ s_H(\gamma') + U-U = s_H(\gamma_0) + U-U$ whereas $t\in \gamma' + \gamma - B_R = \gamma_0 + B_R$. Hence $(w,t)\in (s_H(\gamma_0), \gamma_0) + (U-U)\times B_R$, as desired. This shows the relative density of $\mathcal{G}(s_H)$ so the proof is complete.
\end{proof}

\vspace{0.2cm}





Therefore we get a CPS $(H,\mathsf{\Gamma}, s_H)$ for which the sets $s_H(\mathsf{\Lambda}_i)$ have compact closure $V_i$ in $H$, and moreover such that $(V_i)_{i\in I}$ forms an irredundant family in $H$. It is in fact true that the compact sets $V_i$ are topologically regular in $H$, something which will be deduced in the next step of our proof.\\

\textbf{Step 3: Connection with the CPS $(\mathsf{H},\mathsf{\Gamma}, s_{\mathsf{H}})$.} The following key proposition establish a connection between the CPS $(H,\mathsf{\Gamma}, s_H)$ we are considering with the one provided by Theorem \ref{theo.construction.CPS}:

\vspace{0.2cm}

\begin{prop}\label{prop.universal.CPS} Suppose we have a CPS $(H, \mathsf{\Gamma}, s_H)$ such that the closure $V_i$ of the sets $s_H(\mathsf{\Lambda}_i)$ are compact in $H$ and form an irredundant family $(V_i)_{i\in I}$ in $H$. Then there exists a continuous, onto and open group morphism 
\begin{align*} \theta : \mathsf{H}\longrightarrow H
\end{align*}

with $s_H= \theta \circ s_{\mathsf{H}} $ on the structure group $\mathsf{\Gamma}$. Moreover the family $(V_i)_{i\in I}$ forms a window in $H$ such that $\theta (\mathsf{W}_i)= V_i$ for each $i\in I$.
\end{prop}

\vspace{0.3cm}

\textit{Proof of proposition \ref{prop.universal.CPS}.} The proof of this proposition we be obtained after a series of four lemmas. Considering the group $\mathsf{\Gamma}$ let
\vspace{0.2cm}
\begin{align*} \Xi ^{\mathsf{\Gamma}} := \left\lbrace \Lambda \in \mathbb{X}_{\mathsf{\Lambda}} \, \vert \; Supp(\Lambda )\subseteq \mathsf{\Gamma} \right\rbrace 
\end{align*}

\vspace{0.2cm}

equipped with the combinatorial topology. The lemma given below is an adaptation to the multiple set setting of a result of Schlottmann \cite{Sch}, lemma 4.1 there:


\vspace{0.2cm}

\begin{lem}\label{lem.Sch} Each $\Lambda \in \Xi ^{\mathsf{\Gamma}}$ defines a unique element $w_\Lambda \in H$ through 
\vspace{0.1cm}
\begin{align*} \left\lbrace w_\Lambda \right\rbrace = \bigcap _{i\in I} \bigcap _{\gamma \in \Lambda _i} s_H(\gamma )-V_i
\end{align*}

\end{lem}
\vspace{0.1cm}

\begin{proof} For each $\Lambda \in \Xi ^{\mathsf{\Gamma}}$ denote $V_i(\Lambda):= \overline{s_H(\Lambda _i)}^H$. With respect to this notation the compact sets $V_i$ are equal to $V_i(\mathsf{\Lambda })$. Then given any $\Lambda \in \Xi ^{\mathsf{\Gamma}}$, one has for $w\in H$ an equivalence of conditions
\begin{align*} w\in \bigcap _{i\in I} \bigcap _{\gamma \in \Lambda _i} s_H(\gamma )-V_i \qquad \Longleftrightarrow \qquad  V_i(\Lambda)\subseteq V_i+ w \; \; \forall \, i\in I
\end{align*}
This is true since, as obviously $w\in s_H(\gamma )-V_i$ if and only if $s_H(\gamma ) \in V_i+w$, the left hand condition is equivalent to have $s_H(\Lambda_i)\subseteq V_i+ w$ for each $i\in I$, this latter condition being equivalent by taking closure with the right hand condition above. Let us show that any $\Lambda \in \Xi ^{\Gamma}$ admits an element $w_\Lambda \in H$ where these equivalent conditions hold:\\

From repetitivity of $\mathsf{\Lambda}$ there exists for any $\Lambda \in \Xi ^{\mathsf{\Gamma}}$ a sequence $(\gamma _n)_n$ in $\mathsf{\Gamma}$ such that $\mathsf{\Lambda}+\gamma _n$ converges to $\Lambda$ with respect to the combinatorial topology on $\Xi ^{\mathsf{\Gamma}}$. Thus for each $\gamma \in \Lambda _i$ we eventually get $\gamma \in \mathsf{\Lambda}+\gamma _n$, ensuring that $s_H(\gamma ) \in V_i+s_H(\gamma _n)$ eventually. Picking one of these $\gamma$, we get that the sequence $s_H(\gamma _n)$ lies into each compact set $s_H(\gamma )- V_i$ eventually, and thus accumulates at some element $w_\Lambda\in H$. We can suppose after possibly extracting that $s_H(\gamma _n)$ converges to $w_\Lambda$ in $H$. The latter must satisfy $s_H(\gamma ) \in V_i+w_\Lambda$ for each $\gamma \in \Lambda _i$ and each $i\in I$, as desired.\\

Such element $w_\Lambda$ is unique: from the above we know that the set $V_i(\Lambda)$, being contained into some translate of the compact set $V_i$, is compact in $H$. Then interchanging the roles of $\Lambda $ and $\mathsf{\Lambda}$ in the previous argument shows that there equally exists some $w_\Lambda '\in H$ with $V_i\subseteq V_i(\Lambda)+ w_\Lambda ' \; \; \forall \, i\in I$. This yields $V_i \subseteq V_i(\Lambda) +w_\Lambda ' \subseteq V_i +w_\Lambda + w_\Lambda '$ with $V_i$ compact, which forces $V_i= V_i(\Lambda)+ w_\Lambda+ w_\Lambda ' \; \; \forall \, i\in I$. From the irredundancy assumption on the family $(V_i)_{i\in I}$ one gets $w_\Lambda= -w_\Lambda '$, and this in turns gives that $V_i(\Lambda)= V_i+ w_\Lambda $ for each $ i\in I$. Such an equality, again from irredundancy of the family $(V_i)_{i\in I}$, can be satisfied by at most one element $w_\Lambda$, giving unicity and thus the proof.
\end{proof}

\vspace{0.2cm}

Therefore whenever we are given a CPS as stated in Proposition \ref{prop.universal.CPS} then we get a mapping $\omega : \Xi ^{\mathsf{\Gamma}} \longrightarrow \, H$ with $\omega(\Lambda):= w_\Lambda$ given by the previous lemma. Each $w_\Lambda$ is also obtained as the unique solution to the equation $\overline{ s_H(\Lambda _i)}^{H}= V_i+ w_\Lambda$.

\vspace{0.2cm}

\begin{lem}\label{lem.omega} The mapping $\omega : \Xi ^{\mathsf{\Gamma}} \longrightarrow \, H$ is uniformly continuous, and satisfy the $\mathsf{\Gamma}$-equivariance condition $\omega (\Lambda -\gamma)= \omega (\Lambda) -s_H(\gamma)$ where $\Lambda \in \Xi ^{\mathsf{\Gamma}}$ and $\gamma \in \mathsf{\Gamma}$.
\end{lem}

\vspace{0.2cm}

\begin{proof} The $\mathsf{\Gamma}$-equivariance of $\omega$ is direct: For, if $\Lambda \in \Xi ^{\mathsf{\Gamma}}$ and $\gamma _0 \in \mathsf{\Gamma}$ then
\vspace{0.2cm}
\begin{align*}w_{\Lambda -\gamma _0}  = \left[ \bigcap _{i\in I} \bigcap _{\gamma \in \Lambda _i- \gamma _0} s_H(\gamma )-V_i \right] = \left[ \bigcap _{i\in I} \bigcap _{\gamma ' \in \Lambda _i} s_H(\gamma ')-V_i \right] - s_H(\gamma _0) = w_{\Lambda} -s_H(\gamma _0)
\end{align*}

\vspace{0.2cm}

Let us show the continuity of the mapping $\omega$: Given a neighborhood $U$ of $0$ in $H$, and any $\Lambda \in \Xi ^{\mathsf{\Gamma}} $, since $w_\Lambda$ satisfies the intersection of Lemma \ref{lem.Sch} there exists a radius $R_\Lambda$ such that
\vspace{0.2cm}
\begin{align*} \left[ \bigcap _{i\in I}\bigcap _{\gamma \in \Lambda _i\cap B_{R_\Lambda}}s_H(\gamma) -V_i \right] \, \subseteq \,  w_\Lambda  +U
\end{align*}

\vspace{0.2cm}

Now for any other $\Lambda '\in \Xi ^{\mathsf{\Gamma}} $ coinciding with $\Lambda $ on the ball $B_{R_\Lambda}$ one has in turn 
\vspace{0.2cm}
\begin{align*}w_{\Lambda '} \in \, \left[  \bigcap _{i\in I}\bigcap _{\gamma \in \Lambda ' _i\cap B_{R_\Lambda}}s_H(\gamma) -V_i\right]  \, = \, \left[  \bigcap _{i\in I}\bigcap _{\gamma \in \Lambda  _i\cap B_{R_\Lambda}}s_H(\gamma) -V_i \right]  \,\subseteq \, w_\Lambda +U
\end{align*}

\vspace{0.2cm}

and thus $w_\Lambda  -w_{\Lambda '}\in U$. This shows that $\omega$ is continuous at any $\Lambda \in \Xi ^{\mathsf{\Gamma}} $. To see that $\omega$ is in fact uniformly continuous on $\Xi ^{\mathsf{\Gamma}} $ observe first that this is the case on the compact subset $\Xi$, that is, for any neighborhood $U$ of $0$ in $H$ there is a radius $R_U$ such that whenever $\Lambda , \Lambda '\in \Xi $ coincide on the Euclidean $R_U$-ball then $w_\Lambda  -w_{\Lambda '}\in U$. Suppose then that $\Lambda , \Lambda '\in \Xi ^{\mathsf{\Gamma}}$ coincide on the Euclidean $R_U +R_0$-ball, with $R_0$ some radius of relative density for point patterns in $\mathbb{X}_{\mathsf{\Lambda}}$: Then one can find some element $\gamma \in \mathsf{\Gamma}$ in the Euclidean $R_0$-ball which simultaneously lies into $\Lambda , \Lambda '$, and thus $\Lambda -\gamma, \Lambda '-\gamma$ lie into $\Xi$ and moreover agree on the Euclidean $R_U$-ball. This ensure that $w_{\Lambda -\gamma}  -w_{\Lambda '-\gamma}\in U$, and since the mapping $\omega$ satisfies the equivariance condition $w_{\Lambda -\gamma}= w_\Lambda -s_H(\gamma)$ one finally gets $w_\Lambda  -w_{\Lambda '}\in U$. This shows the uniform continuity of $\omega$, concluding the proof.
\end{proof}

\vspace{0.2cm}








\begin{lem}\label{lem.factorization} The mapping $\omega$ factorizes into $\displaystyle \Xi ^{\mathsf{\Gamma}} \longrightarrow^{[.]_{srp}} \, \mathsf{H} \longrightarrow ^{\theta} \, H$ where $[.]_{srp}$ is the quotient map under the strong regional proximality relation and $\theta$ is a continuous, onto and open group morphism with $s_{H}= \theta \circ s_{\mathsf{H} }$.
\end{lem}

\vspace{0.2cm}

\begin{proof} Recall that $[\Lambda]_{srp}\subseteq \Xi ^{\mathsf{\Gamma}}$ whenever $\Lambda \in \Xi ^{\mathsf{\Gamma}}$ by Proposition \ref{prop.construction.CPS}. The uniform continuity of the mapping $\omega$ shown in the previous lemma ensure that whenever $[\Lambda]_{srp} = [\Lambda ']_{srp}$ then $w_\Lambda  =w_{\Lambda '}$: For, if $[\Lambda]_{srp} = [\Lambda ']_{srp}$ then for each radius $R$ one can find $\Lambda _1, \Lambda _2$ and $t\in \mathbb{R}^d$, this latter being possibly supposed to lie into the relatively dense subgroup $\mathsf{\Gamma}$, such that 
\begin{align*} \Lambda \cap B_R &=  \Lambda _1 \cap B_R \\ \Lambda ' \cap B_R &=  \Lambda _2 \cap B_R \\ (\Lambda _1-t)\cap B_R &=  (\Lambda _2-t) \cap B_R
\end{align*}

\vspace{0.2cm}
One easily show that the point patterns  $\Lambda _1$ and $ \Lambda _2$ must be supported on $\mathsf{\Gamma}$, so that $\Lambda _1, \Lambda _2\in \Xi ^{\mathsf{\Gamma}}$ and thus $w_{\Lambda _1}$ and  $w_{\Lambda _2}$ are well defined elements in $H$. For each neighborhood $U$ of $0$ in $H$ there is an $R$ such that the three equalities above give
\begin{align*} w_{\Lambda } - w_{\Lambda _1} &\in U \\ w_{\Lambda '} - w_{\Lambda _2} &\in U \\ w_{\Lambda _1} - w_{\Lambda _2} &= w_{\Lambda _1-t} - w_{\Lambda _2-t}\in U
\end{align*}

\vspace{0.2cm}

which yields $w_{\Lambda }-w_{\Lambda '}\in U-U+U$ for any $U$. This forces $w_{\Lambda }=w_{\Lambda '}$, as desired.
\\



Hence the mapping $\omega$ factorizes into $\displaystyle \Xi ^{\mathsf{\Gamma}} \longrightarrow^{[.]_{srp}} \, \mathsf{H} \longrightarrow ^{\theta} \, H$ where $[.]_{srp}$ is the quotient map under the strong regional proximality relation. By definition of the quotient topology on $\mathsf{H}$ the mapping $\theta$ of this factorization is continuous. Now by definition of $s_{\mathsf{H}}$ one has $\theta (s_{\mathsf{H}}(\gamma))= \theta ([\mathsf{\Lambda} +\gamma]_{srp})= \omega (\mathsf{\Lambda} +\gamma)= \omega(\mathsf{\Lambda})+s_{H}(\gamma)= s_{H}(\gamma)$, which yields $s_{H}= \theta \circ s_{\mathsf{H}}$. As both $\mathsf{H}$ and $H$ are LCA group completions of the group $\mathsf{\Gamma}$ the mapping $\theta$ must in turn be a group morphism. Finally, $\theta$ must be onto and open from lemma \ref{lem.CPS}.
\end{proof}

\vspace{0.2cm}





\begin{lem} Each $V_i$ is topologically regular in $H$, with $\theta (\mathsf{W}_i)= V_i$ for each $i\in I$.

\end{lem} 

\vspace{0.2cm}

\begin{proof} Let us show that $\theta (\mathsf{W}_i)= V_i$ for each $i\in I$: From point $(iii)$ of Theorem \ref{theo.construction.CPS} we have $\mathsf{W}_i= - [\Xi _i]_{srp}$ in $\mathsf{H}$, and consequently $\theta (\mathsf{W}_i)= - \omega(\Xi _i)$. By repetitivity of $\mathsf{\Lambda}$ the set $\left\lbrace \mathsf{\Lambda} -\gamma \, \vert \gamma \in \mathsf{\Lambda}_i \right\rbrace $ is dense in the transversal $\Xi _i$, and because $\omega$ is continuous one gets that $ \omega(\Xi _i)$ is equal to the closure in $H$ of $\omega(\left\lbrace \mathsf{\Lambda} -\gamma \, \vert \gamma \in \mathsf{\Lambda}_i \right\rbrace )$. However $\omega( \mathsf{\Lambda} -\gamma)= -s_H(\gamma)$ by lemma \ref{lem.omega} so the latter is nothing but the closure in $H$ of $- s_H(\mathsf{\Lambda}_i)$, which is exactly $-V_i$. This gives $ V_i = - \omega(\Xi _i) = \theta (\mathsf{W}_i)$, as desired. Now using this, together with the fact that $\theta $ is an open map from previous lemma, we obtain $\theta (\mathring{\mathsf{W}}_i)\subseteq \mathring{V}_i$ for each $i\in I$. Since $\mathsf{W}_i$ is the closure of its interior then $\theta (\mathring{\mathsf{W}}_i)$ must be dense in $\theta (\mathsf{W}_i)= V_i$. This ensure that each $V_i$ is the closure of its interior, hence topologically regular in $H$.
\end{proof}

\vspace{0.2cm}

These four lemmas settle the proof of Proposition \ref{prop.universal.CPS}.
\begin{flushright}
$\square $
\end{flushright}

\vspace{0.2cm}


\textbf{Step 4: Concluding argument.} We shall now connect the model multiple set $\underline{\mathsf{\Lambda}}$ stated at the beginning of this proof with $\Delta$:\\

First one has for each $i\in I$ the inclusion of $\mathfrak{P}_{\mathsf{H}}(\mathring{\mathsf{W}}_i)$ into $\mathfrak{P}_{H}(\mathring{V}_i)$: For, as the morphism $\theta : \mathsf{H} \longrightarrow H$ provided by proposition \ref{prop.universal.CPS} is open and maps $\mathsf{W}_i$ onto $V_i$ one has $\theta (\mathring{\mathsf{W}}_i)\subseteq \mathring{V}_i$. Now has $s_H = \theta \circ s_{\mathsf{H}}$ we get that 
\vspace{0.2cm}
\begin{align*}\gamma \in \mathfrak{P}_{\mathsf{H}}(\mathring{\mathsf{W}}_i) \; \Longleftrightarrow \; s_{\mathsf{H}}(\gamma) \in \mathring{\mathsf{W}}_i \; \Longrightarrow \; s_H(\gamma) \in \mathring{V}_i \; \Longleftrightarrow \; \gamma \in \mathfrak{P}_{H}(\mathring{V}_i)
\end{align*}

\vspace{0.2cm}

Next we shall connect the model multiple set $(\mathfrak{P}_{H}(\mathring{V}_i))_{i\in I}$ with CPS $(H^2, \Gamma^2, s_{H^2})$ and collection $(V_i^2)_{i\in I}$: The couple of CPS $(H, \mathsf{\Gamma}, s_H)$ and $(H^2, \Gamma^2, s_{H^2})$ is so that $H \subseteq H^2$, $\mathsf{\Gamma} \subseteq \Gamma ^2$, and the $*$-map $s_H$ is the restriction of $s_{H_2}$ on $\mathsf{\Gamma}$, and thus one can apply Proposition \ref{prop.general.CPS} to get that the interior of each $V_i$ in $H$ is equal to the interior of $V_i^2$ relatively to the larger group $H^2$. This in turn yields that $ \mathfrak{P}_{H}(\mathring{V}_i) \subseteq \mathfrak{P}_{H^2}(\mathring{V}_i^2) $.\\

Now the CPS $(H^2, \Gamma^2, s_{H^2})$ with collection $(V_i^2)_{i\in i}$ is connected with the CPS $(H^1, \Gamma^1, s_{H^1})$ with collection $(V_i^1)_{i\in I}$ by the group morphism $[.]_\mathfrak{R}: H^2 \longrightarrow H^1$ modding out the redundancy subgroup $\mathfrak{R}$ of the family $(V_i^1)_{i\in I}$ by $s_{H^2} = [.]_\mathfrak{R} \circ s_{H^1}$ (see formula (\ref{morphism2})) and $[.]_\mathfrak{R}^{-1}(V_i^2)= V_i^1$. It follows that 
\vspace{0.2cm}
\begin{align*}\gamma \in \mathfrak{P}_{H^2}(\mathring{V}_i^2) \; \Longleftrightarrow \; s_{H^2}(\gamma) \in \mathring{V}_i^2 \; \Longrightarrow \; s_{H^1}(\gamma) \in [.]_\mathfrak{R}^{-1}(\mathring{V}_i^2)\subseteq  \mathring{V}_i^1 \; \Longleftrightarrow \; \gamma \in \mathfrak{P}_{H^1}(\mathring{V}_i^1)
\end{align*}

\vspace{0.2cm}

Now each set $V_i^1$ is included into the compact set $W_i^1$ and thus $\mathfrak{P}_{H^1}(\mathring{V}_i^1)\subseteq \mathfrak{P}_{H^1}(\mathring{W}_i^1)$ for each $i\in I$. Above all, one gets
\vspace{0.2cm}
\begin{align*}\mathsf{\Lambda}_i \; \subseteq \; \underline{\mathsf{\Lambda}}_i &  = \; \mathfrak{P}_{\mathsf{H}}(\mathring{\mathsf{W}}_i)\cup \mathsf{\Lambda}_i \\ &  \subseteq \; \mathfrak{P}_{H}(\mathring{V}_i)\cup \mathsf{\Lambda}_i\\  &  \subseteq  \; \mathfrak{P}_{H^2}(\mathring{V}_i^2)\cup \mathsf{\Lambda}_i \\  & \subseteq  \; \mathfrak{P}_{H^1}(\mathring{V}_i^1) \cup \mathsf{\Lambda}_i \\  & \subseteq \; \mathfrak{P}_{H^1}(\mathring{W}_i^1)\cup \mathsf{\Lambda}_i  \subseteq \Delta_i
\end{align*}

\vspace{0.2cm}

whenever $i\in I$, concluding the proof of Theorem \ref{theo.principal}.
\begin{flushright}
$\square$
\end{flushright}

\section{On morphisms between CPS}

In the process of the proof of Theorem \ref{theo.principal} it has been used general facts about CPS, which we prove in this dedicated section. Let us first recall a general fact:

\vspace{0.3cm}

\begin{theo}\label{theo.open.map}(Open Mapping Theorem for locally compact groups) Let $\theta : G^1 \longrightarrow G^2$ be a group morphism between locally compact groups. If $\theta$ is surjective and $G^1$ is $\sigma$-compact, then $\theta$ is an open map.
\end{theo}

\vspace{0.3cm}
The general theorem above will serve to show two openness properties for CPS, namely the lemma just below which is invoked in the proof of Lemma \ref{lem.factorization}, and a proposition given next which is invoked at the final step of our proof of Theorem \ref{theo.principal}.

\vspace{0.2cm}

\begin{lem}\label{lem.CPS} Let $(H_1,\Gamma, s_1)$ and $(H_2,\Gamma, s_2)$ be two CPS with same structure group $\Gamma$. Suppose we have a continuous group morphism $\theta : H_1 \longrightarrow H_2$ such that $s_2= \theta \circ s_1$. Then $\theta$ is onto and open.

\end{lem}

\vspace{0.2cm}

\begin{proof} The continuous morphism $\theta$ provides a commutative diagram
\vspace{0.2cm}
\begin{align*}\begin{psmatrix}[colsep=1.5cm,
rowsep=1.2cm]
H_1\times \mathbb{R}^d \; \; & \; \; H_2\times \mathbb{R}^d\\
\left[ H_1\times \mathbb{R}^d\right] _\Gamma\; \;   & \; \; \; \left[ H_2\times \mathbb{R}^d\right] _\Gamma
\psset{arrows=->>,linewidth=0.2pt, labelsep=1.5pt
,nodesep=0pt}
\ncline{1,1}{2,1}<{[.]_{\Gamma}}
\ncline{1,2}{2,2}>{[.]_{\Gamma}}
\ncline{2,1}{2,2}^{\tilde{\theta}}
\psset{arrows=->,linewidth=0.2pt, labelsep=1.5pt
,nodesep=0pt}
\ncline{1,1}{1,2}^{\theta \times id}
\end{psmatrix}\end{align*}

\vspace{0.3cm}
where vertical arrows are quotient morphisms under the graphs subgroups of the $*$-maps $s_1$ and $s_2$ respectively, and the lower arrow is a continuous morphism given by $\tilde{\theta}([w,t]_\Gamma):= [\theta(w),t]_\Gamma$. Now $\tilde{\theta}$ must be onto: indeed its domain $\left[ H_1\times \mathbb{R}^d\right] _\Gamma$ is a compact space and has a range in $\left[ H_2\times \mathbb{R}^d\right] _\Gamma$ containing at least $\left[ \left\lbrace 0\right\rbrace \times \mathbb{R}^d\right] _\Gamma$. But this latter is dense since the $*$-map $s_2$ has by assumption a dense range $s_2(\Gamma)$ in $H_2$. Therefore the map $\tilde{\theta}$ has a range which is either compact and dense in $\left[ H_2\times \mathbb{R}^d\right] _\Gamma$, and thus is onto.\\

Now we prove the surjectivity of $\theta$: If $w\in H_2$ is given then one has $[w,0]_\Gamma \in \left[ H_2\times \mathbb{R}^d\right] _\Gamma$ which admits a antecedent class under $\tilde{\theta}$, that is, there is $\tilde{w}\in H_1$ and $t\in \mathbb{R}^d$ such that $\tilde{\theta}([\tilde{w},t]_\Gamma) = [w,0]_\Gamma$. Hence there exists some $\gamma\in \Gamma$ such that $(\theta(\tilde{w})+s_2(\gamma),t+\gamma)= (w,0)$ in $H_2\times \mathbb{R}^d$. Thus the element $w':=\tilde{w}+s_1(\gamma)\in H_1$ is so that $\theta (w')= \theta(\tilde{w})+s_2(\gamma)=w$, proving the surjectivity. The openess of $\theta$ follows from the open mapping theorem (Theorem \ref{theo.open.map}), finishing the proof.
\end{proof}

\vspace{0.2cm}

\begin{prop}\label{prop.general.CPS} Let $(H_1,\Gamma_1, s_{H_1})$ and $(H_2,\Gamma_2, s_{H_2})$ be two CPS such that $\Gamma_1\subseteq \Gamma_2$, and that there is a continuous group morphism $\theta : H_1 \longrightarrow H_2$ with $s_{H_2}= \theta \circ s_{H_1}$ on $\Gamma_1$. Then $\theta$ is open and $\theta (H_1)$ is an open subgroup of $H_2$.
\end{prop}

\vspace{0.2cm}

\begin{proof} Consider $\Gamma_2$ endowed with the discrete topology, and let $H_1\oplus _{\, \Gamma _1} \Gamma_2$ be the quotient of the LCA group $H_1\oplus \Gamma_2$ by the discrete subgroup $\left\lbrace (-s_{H_1}(\gamma), \gamma) \, \vert \, \gamma \in \Gamma_1\right\rbrace $. There is a well defined group morphism
\vspace{0.2cm}
\begin{align}\label{morphism}\begin{psmatrix}[colsep=1.5cm,
rowsep=1.2cm]
\Gamma_2 \ni \gamma \; \; & \; \; \; [0,\gamma]_{\Gamma_1} \in H_1\oplus _{\, \Gamma _1} \Gamma_2
\psset{arrows=|->,linewidth=0.2pt, labelsep=2pt
,nodesep=0pt}
\ncline{1,1}{1,2}^{s_{H_1\oplus _{\, \Gamma _1} \Gamma_2}\qquad}
\end{psmatrix}\end{align}

\vspace{0.2cm}

Observe that in case $\Gamma_1= \Gamma_2$ then $H_1\oplus _{\, \Gamma _1} \Gamma_2= H_1$ and $s_{H_1\oplus _{\, \Gamma _1} \Gamma_2}= s_{H_1}$. In the general case, one has a continuous group morphism 
\vspace{0.2cm}
\begin{align*}\Theta : H_1\oplus _{\, \Gamma _1} \Gamma_2 \longrightarrow H_2 \qquad \qquad \Theta ([w,\gamma]_{\Gamma_1}):= \theta(w)+s_{H_2}(\gamma)
\end{align*}

\vspace{0.2cm}

satisfying $s_{H_2}= \Theta \circ s_{H_1\oplus _{\, \Gamma _1} \Gamma_2}$ on $\Gamma_2$. Now the triple $(H_1\oplus _{\, \Gamma _1} \Gamma_2, \Gamma_2, \mathbb{R}^d)$ obtained is in fact a CPS with $*$-map $s_{H_1\oplus _{\, \Gamma _1} \Gamma_2}$. For, first of all the group morphism $s_{H_1\oplus _{\, \Gamma _1} \Gamma_2}$ has dense range since given $(w,\gamma)\in H_1 \oplus \Gamma_2$ one has a sequence $(\gamma_n)_n\subset \Gamma_1$ such that $s_{H_1}(\gamma_n)$ converges to $w$ in $H_1$, and thus $s_{H_1\oplus _{\, \Gamma _1} \Gamma_2}(\gamma-\gamma_n) = [0,\gamma -\gamma_n]_{\Gamma_1}= [s_H(\gamma_n),\gamma]_{\Gamma_1}$ converges to $[w,\gamma]_{\Gamma_1}$ in $H_1\oplus _{\, \Gamma _1} \Gamma_2$. Second the graph $\mathcal{G}(s_{H_1\oplus _{\, \Gamma _1} \Gamma_2})$ is uniformly discrete in $H_1\oplus _{\, \Gamma _1} \Gamma_2\times \mathbb{R}^d$ since one has an obvious continuous group morphism
\vspace{0.2cm}
\begin{align*}\Theta \times id : H_1\oplus _{\, \Gamma _1} \Gamma_2\times \mathbb{R}^d \longrightarrow H_2\times \mathbb{R}^d
\end{align*}

\vspace{0.2cm}

mapping $\mathcal{G}(s_{H_1\oplus _{\, \Gamma _1} \Gamma_2})$ onto the uniformly discrete subgroup $\mathcal{G}(s_{H_2})$ of $H_2\times \mathbb{R}^d$. Third $\mathcal{G}(s_{H_1\oplus _{\, \Gamma _1} \Gamma_2})$ is relatively dense in $H_1\oplus _{\, \Gamma _1} \Gamma_2\times \mathbb{R}^d$ since whenever $K\subset H_1$ and $K'\subset \mathbb{R}^d$ are compact subsets with $\mathcal{G}(s_{H_1})+K\times K' = H_1\times \mathbb{R}^d$ then $[K,0]_{\Gamma_1}\times K'$ is compact with 
\vspace{0.2cm}
\begin{align*}\mathcal{G}(s_{H_1\oplus _{\, \Gamma _1} \Gamma_2})+ [K,0]_{\Gamma_1}\times K'= H_1\oplus _{\, \Gamma _1} \Gamma_2\times \mathbb{R}^d
\end{align*}

\vspace{0.2cm}

For, if $(w, \gamma, t)$ is given then there is $\gamma ' \in  \Gamma_1 $ with $(w,t-\gamma)\in (K+s_H(\gamma'), K'+\gamma')$, and thus considering $\gamma_0:= \gamma + \gamma' \in \Gamma^2$ provides on one hand
\vspace{0.2cm}
\begin{align*} [w,\gamma]_{\Gamma_1} = [w-s_H(\gamma'),\gamma +\gamma']_{ \Gamma_1}= [w-s_H(\gamma'),\gamma_0]_{\Gamma_1}\in [K,\gamma_0]_{\Gamma_1}= [0,\gamma_0]_{ \Gamma _1}+ [K,0]_{\Gamma_1}
\end{align*}

\vspace{0.2cm}

and on the other hand gives $t\in \gamma + \gamma' +K' = \gamma_0 + K'$, yielding
\vspace{0.2cm}
\begin{align*} ([w,\gamma]_{\Gamma_1},t)\in ([0,\gamma_0]_{ \Gamma_1}, \gamma_0)+ [K,0]_{\Gamma_1}\times K'
\end{align*}

\vspace{0.2cm}

showing the relative density of $\mathcal{G}(s_{H_1\oplus _{\, \Gamma _1} \Gamma_2})$ in $H_1\oplus _{\, \Gamma _1} \Gamma_2\times \mathbb{R}^d$. Therefore we have a new CPS $(H_1\oplus _{\, \Gamma _1} \Gamma_2, \Gamma_2, \mathbb{R}^d)$ with a continous group morphism $\Theta : H_1\oplus _{\, \Gamma _1} \Gamma_2 \longrightarrow H_2$ with  $s_{H_2}= \Theta \circ s_{H_1\oplus _{\, \Gamma _1} \Gamma_2}$ on $\Gamma_2$. Lemma \ref{lem.CPS} guarantee that $\Theta$ is an open map, that is, the image of any open subset of $H_1\oplus _{\, \Gamma _1} \Gamma_2$ is open in $H_2$. However the copy $[H_1,0]_{\Gamma_1} $ of $H_1$ is open in $H_1\oplus _{\, \Gamma _1} \Gamma_2$: For it is the image of the open subgroup $H_1\times \left\lbrace 0\right\rbrace  $ of $H_1\oplus \Gamma_2$ under the onto, hence open, continuous group morphism taking the quotient by the discrete subgroup $\left\lbrace (-s_{H_1}(\gamma), \gamma) \, \vert \, \gamma \in \Gamma_1\right\rbrace $. We conclude by observing that the formula defining $\Theta$ gives $\theta (H_1)= \Theta ([H_1,0]_{\Gamma_1} )$ and so is open in $H_2$, and that the image under $\theta$ of any open subset $U$ of $H_1$ is equal to the image under $\Theta$ of the open subset $[U,0]_{\Gamma_1}$ of $H_1\oplus _{\, \Gamma _1} \Gamma_2$ and thus is open in $H_2$.
\end{proof}

\vspace{0.3cm}

\textbf{Acknowledgments.} The author thanks the Universidad de Santiago de Chile, USACH and Proyecto POSTDOC$\_$DICYT C\'odigo 001316POSTDOC, Vicerrectoria de Investigaci\'on, Desarrollo e Innovaci\'on.

\end{document}